\documentclass[11pt]{article}
\usepackage[T1]{fontenc}
\usepackage[utf8]{inputenc}
\usepackage{lmodern}
\usepackage{microtype}
\usepackage{amsmath,amssymb,amsthm,mathtools}
\usepackage[margin=2.4cm]{geometry}
\usepackage{booktabs}
\usepackage{parskip}           % no paragraph indentation
\usepackage[backref=page]{hyperref}
\hypersetup{colorlinks=true,linkcolor=blue,urlcolor=blue}

\hypersetup{hidelinks}

\newtheorem{theorem}{Theorem}[section]
\numberwithin{equation}{section}

\newtheorem{lemma}[theorem]{Lemma}
\newtheorem{corollary}[theorem]{Corollary}
\newtheorem{introthm}{Theorem}            % main results in the Introduction: A, B, C, D, E

\theoremstyle{definition}
\newtheorem{definition}[theorem]{Definition}
\newtheorem{remark}[theorem]{Remark}
\newtheorem{example}[theorem]{Example}

\newcommand{\BR}{\mathcal{BR}}
\newcommand{\YB}{\mathcal{YB}}

\usepackage{enumerate}
\newcommand{\selabel}[1]{\label{se:#1}}
\newcommand{\seref}[1]{Section~\ref{se:#1}}

\def\thlabel#1{\label{th:#1}}
\def\thref#1{Theorem~\ref{th:#1}}

\newcommand{\lelabel}[1]{\label{le:#1}}
\newcommand{\leref}[1]{Lemma~\ref{le:#1}}

\newcommand{\colabel}[1]{\label{co:#1}}
\newcommand{\coref}[1]{Corollary~\ref{co:#1}}
\newcommand{\relabel}[1]{\label{re:#1}}
\newcommand{\reref}[1]{Remark~\ref{re:#1}}
\newcommand{\exlabel}[1]{\label{ex:#1}}

\newcommand{\delabel}[1]{\label{de:#1}}

\newcommand{\eqlabel}[1]{\label{eq:#1}}
\newcommand{\equref}[1]{(\ref{eq:#1})}

\newtheorem{remarks}[theorem]{Remarks}

\newcommand{\bZ}{\mathbb{Z}}

\title{Solutions of the set-theoretic quantum Yang--Baxter equation arising from posets}
\author{G. Militaru}
\date{}

\begin{document}
\maketitle

\begin{abstract}
Victoria Lebed observed that the map $(x, \, y) \mapsto (\min\{x, \, y\}, \, \max\{x, \, y\})$ on a totally ordered set is an idempotent solution of the braid equation. We extend this example in several directions, by constructing four families of solutions of the braid equation associated to posets. First, for a given poset $P$, the map $R_P$ which sorts the comparable pairs and fixes the incomparable ones is an idempotent solution of the braid equation whose automorphism group is the automorphism group of the poset $P$. Combined with the classical theorem of Birkhoff, this shows that every group is isomorphic to the automorphism group of an idempotent solution of the braid equation of the form $(P, \, R_P)$, for some poset $P$. Moreover, $R_P \cong R_Q$ if and only if $P \cong Q$, and the solutions of the form $R_P$ are exactly the idempotent solutions of the braid equation which permute the components of each pair. Secondly, we classify all the solutions $R$ of the braid equation on an arbitrary set $X$ such that $R(x, \, y) \in \{(x, \, y), \, (y, \, x)\}$
for all $x$, $y\in X$: they are exactly the solutions $R_{\mathcal{P}}$ associated to the decompositions $\mathcal{P}$ of $X$ into posets. Finally, we associate to any finite poset $P$ two involutive solutions $D_P$ and $E_P$ of the braid equation, whose automorphism groups are $(\mathbb{Z}/2)^P \rtimes {\rm Aut} (P)$ and, if $P$ is connected and has at least two elements, the group of automorphisms and anti-automorphisms of $P$.
\end{abstract}

\medskip

\noindent {\em Key words: set-theoretic Yang--Baxter equation; braid equation; poset; idempotent solution; involutive solution; automorphism group}

\vspace{.3cm}

\noindent{MSC 2020: Primary 16T25; Secondary 06A06, 20B25, 20B27, 08A35}

\section*{Introduction}

Few equations have connected so many different areas of mathematics as the Yang--Baxter equation. Introduced in statistical mechanics and later central to the theory of Hopf algebras and quantum groups \cite{K}, it was reduced by Drinfeld \cite{Drinfeld92} to its simplest, purely combinatorial, form: find all the maps $r \colon X \times X \to X \times X$, where $X$ is a set, such that
\begin{equation*}
r^{12} r^{13} r^{23} = r^{23} r^{13} r^{12},
\end{equation*}
or, equivalently, all the maps $R = \tau_X \circ r$ which satisfy the braid equation $R^{12} R^{23} R^{12} = R^{23} R^{12} R^{23}$. Although it looks elementary, this problem has turned out to be related to semigroups \cite{acm1, CMT, GIVdB98}, to invariants of knots and links \cite{BES, CES04}, to the structure groups of solutions \cite{ESS99, LYZ00}, to bijective $1$-cocycles and their equivalent description in terms of braces and skew braces \cite{CV, Rump07}, to combinatorics and graph theory \cite{GIM, MilAut, MilBIP, Yang16}, and to other varieties of algebras, in the sense of universal algebra \cite{ACM2, ACMv3, CMv1}, to mention only a few of the areas in which it appears today. For further details and additional motivation, we refer to the recent monographs \cite{BES, CV} and the references therein.

Most of the literature deals with non-degenerate solutions, but in the last years degenerate solutions, and in particular idempotent ones, have also attracted attention \cite{Lebed17, SV21}. Order structures enter the picture in a natural way: Lebed \cite[Examples 2.7 and 2.8]{Lebed17} obtained idempotent solutions from totally ordered sets and from distributive lattices, and skew lattices produce large classes of solutions \cite{CVV20}.

In this paper we show that arbitrary posets are a natural source of solutions of the braid equation, and that these solutions retain a large part of the information contained in the order. Our work is part of a broader program, announced in \cite[Problems A and B, page 4]{CMv1} (see also \cite{ACM2, acm1}), which can be summarized as follows: first, to identify other categories $\mathcal{C}$ from which YB-solutions can be constructed, ideally in a functorial way; secondly, and this is the more difficult part, to describe precisely the subcategory of YB-solutions obtained by such a construction, i.e. the essential image of the corresponding functor. We believe that this is a natural way to approach Drinfeld's problem. Indeed, the recent classification results show that, at least on sets with a small number of elements, the number of isomorphism classes of solutions arising from left skew braces, which have been the main tool in the study of the Yang--Baxter equation during the last decade, is very small compared to the total number of solutions, and even to the number of non-degenerate ones (see \cite[Section 3]{MilBIP}).

Recently, several families of YB-solutions have been constructed from undirected graphs, with the aim of proving that every group is isomorphic to the automorphism group of a set-theoretic solution of the Yang--Baxter equation \cite{MilAut}. In this paper we use posets instead of graphs. In this case the solutions live on the poset itself, and they are defined in the simplest possible way: by sorting. Our starting point is an example of Lebed \cite[Example 2.7]{Lebed17}: if $X$ is a totally ordered set, then the map $(x, \, y) \mapsto (\min\{x, \, y\}, \, \max\{x, \, y\})$ is an idempotent solution of the braid equation. This example can be extended to arbitrary posets: for a poset $P$, let $R_P \colon P \times P \to P \times P$ be the map which sorts the pairs of comparable elements and fixes the pairs of incomparable elements (see \equref{RP}). The first main result of the paper is the following (\thref{RP} and \coref{caractRP}):

\begin{introthm}
Let $P$ and $Q$ be two posets. Then $R_P$ is an idempotent solution of the braid equation, ${\rm Aut} (P, \, R_P) = {\rm Aut} (P)$, and $(P, \, R_P) \cong (Q, \, R_Q)$ if and only if $P \cong Q$. Conversely, if $X$ is a set and $R \colon X \times X \to X \times X$ is a map such that $R (x, \, y) \in \{(x, \, y), \, (y, \, x)\}$, for all $x$, $y\in X$, then $R$ is an idempotent solution of the braid equation if and only if $R = R_P$, for a unique partial order on $X$.
\end{introthm}

Combined with the theorem of Birkhoff \cite[Teorema 2]{Birkhoff46}, Theorem A gives another proof of the main result of \cite{MilAut}: any group is isomorphic to the automorphism group of an idempotent solution of the braid equation of the form $R_P$, which can be chosen finite if the group is finite (\coref{poset}). Moreover, the classification of the idempotent solutions of the braid equation which permute the components of each pair is equivalent to the classification of posets; in particular, already counting them is a hard problem, for which no closed formula is known (\reref{counting}). Since they permute the components of each pair, both the solutions $R_P$ and their generalizations in Theorem B below are solutions of Lyubashenko type, in the sense of \cite[Definition 2.10]{ACM2}.

The category $\mathcal{P}os$ of posets and strictly monotone maps is alg-universal \cite[Chapter IV, 5.1 and Theorem 5.6]{PultrTrnkova80}; hence the category $\YB$ of all YB-solutions, which is a variety of algebras \cite{CMv1} and in particular an algebraic category in the sense of \cite[Chapter I, 8.4]{PultrTrnkova80}, admits a full embedding, i.e. a full and faithful functor which is injective on objects,
\begin{equation}\eqlabel{YBinPos}
\iota \colon \YB \cong \BR \hookrightarrow \mathcal{P}os,
\end{equation}
where $\BR$ is the category of all solutions of the braid equation. Theorem A goes in the opposite direction, in a completely explicit way: $P \mapsto (P, \, R_P)$ is a full embedding
\begin{equation}\eqlabel{PosinYB}
\Phi \colon \mathcal{P}os^{\bullet} \hookrightarrow \BR \cong \YB
\end{equation}
into $\BR$ of the category $\mathcal{P}os^{\bullet}$ of posets whose morphisms $P \to Q$ are the maps $\sigma$ satisfying
\begin{equation*}
y < x \,\, \Longrightarrow \,\, \sigma (y) \leq \sigma (x) \qquad {\rm and} \qquad \sigma (y) < \sigma (x) \,\, \Longrightarrow \,\, y < x,
\end{equation*}
for all $x$, $y\in P$ (\thref{RP}(b)).\footnote{The category $\mathcal{P}os^{\bullet}$ is different from $\mathcal{P}os$: any constant map satisfies these conditions, while a strictly monotone map need not satisfy the second one.} Thus, for posets, the construction $P \mapsto (P, \, R_P)$ gives an answer to both questions of the program described above: it is a full embedding of the category $\mathcal{P}os^{\bullet}$ into $\mathcal{BR} \cong \YB$, and its essential image is described intrinsically by \coref{caractRP}. Of course, this is only one of the possible ways of constructing YB-solutions from posets; two other constructions are given in \seref{invol}.

The second main result (\thref{flipfix}) is the classification of all the solutions $R$ of the braid equation which permute the components of each pair, i.e. $R(x, \, y) \in \{(x, \, y), \, (y, \, x)\}$ for all $x$, $y$: no other assumption on $R$, such as idempotency, bijectivity or finiteness of the underlying set, is needed.

\begin{introthm}
Let $X$ be a set and $R \colon X \times X \to X \times X$ a map such that $R (x, \, y) \in \{(x, \, y), \, (y, \, x)\}$, for all $x$, $y\in X$. Then $R$ is a solution of the braid equation if and only if there exist a partition $X = \bigsqcup_{i\in I} X_i$ and a partial order on each $X_i$ such that $R$ sorts the comparable pairs of elements of each $X_i$, fixes the incomparable ones, and $R (x, \, y) = (y, \, x)$ whenever $x$ and $y$ belong to different blocks. Moreover, $R$ is bijective if and only if all the blocks $X_i$ are antichains.
\end{introthm}

In particular, a bijective solution of this type depends only on a partition of $X$. In order to obtain bijective solutions which take into account the order of a poset $P$, we have to enlarge
the underlying set; in \seref{invol} we use for this purpose the bipartite solutions of \cite{MilBIP} and we prove (\thref{autDP} and \thref{autEP}):

\begin{introthm}
Let $P$ be a finite non-empty poset and $m$ the number of pairs $(y, \, x) \in P \times P$ with $y < x$. Then there exist two involutive solutions $D_P$ and $E_P$ of the braid equation, on sets $X_P$ and $X'_P$ with $3 |P|$ and $3 |P| + m$ elements, respectively, such that ${\rm Aut} (X_P, \, D_P) \cong (\bZ/2)^P \rtimes {\rm Aut} (P)$ and, if $P$ is connected and has at least two elements, ${\rm Aut} (X'_P, \, E_P)$ is isomorphic to the group of all automorphisms and anti-automorphisms of $P$. In particular, if $P$ is connected and $P \not\cong P^{\rm op}$, then ${\rm Aut} (X'_P, \, E_P) \cong {\rm Aut} (P)$.
\end{introthm}

\section{Preliminaries}\selabel{prel}

Throughout, for a set $X$ we write ${\rm Sym} (X)$ for the group of all permutations of $X$ and $\tau_X \colon X \times X \to X \times X$, $\tau_X (x, \, y) := (y, \, x)$, for the flip map. Let
$R \colon X \times X \to X \times X$ be a map. On $X^3 = X \times X \times X$ we consider the maps
\begin{equation*}
R^{12} := R \times {\rm Id}_X, \qquad R^{23} := {\rm Id}_X \times R, \qquad R^{13} := ({\rm Id}_X \times \tau_X) \circ R^{12} \circ ({\rm Id}_X \times \tau_X).
\end{equation*}
Then $R$ is called a \emph{YB-solution} (i.e. a set-theoretic solution of the \emph{Yang--Baxter equation}), respectively a \emph{solution of the braid equation}, if
\begin{equation}\label{eq:qybe}
R^{12} R^{13} R^{23} = R^{23} R^{13} R^{12}, \qquad \text{respectively} \qquad R^{12} R^{23} R^{12} = R^{23} R^{12} R^{23}.
\end{equation}
It is well known that the braid equation holds for $R$ exactly when $\tau_X \circ R$ satisfies the Yang--Baxter equation. A map $R \colon X\times X \to X\times X$ is called idempotent if $R^2 = R$, involutive if $R^2 = {\rm Id}_{X \times X}$, and unitary if $R^{21} R = {\rm Id}_{X \times X}$, where $R^{21} := \tau_X R \tau_X$. It is often convenient to describe $R$ by its two components. Writing
\begin{equation}\eqlabel{notatieR}
R(x, \, y) = (x \cdot y, \, x \ast y),
\end{equation}
we obtain two binary operations on $X$; in this way, the maps $R \colon X \times X \to X \times X$ are in bijection with the bi-magma structures $(X, \, \cdot, \, \ast)$ on $X$. For $x \in X$, the maps $l_x$, $r_x \colon X \to X$, $l_x (y) := x \cdot y$ and $r_x (y) := y \ast x$, are the left and the right translations of this bi-magma, so that $R (x, \, y) = \bigl( l_x (y), \, r_y (x) \bigr)$.

If $R$ and $R'$ are maps as above on the sets $X$ and $X'$, a morphism from $(X, \, R)$ to $(X', \, R')$ is a map $\sigma \colon X \to X'$ which satisfies $(\sigma \times \sigma) \circ R = R' \circ (\sigma \times \sigma)$; equivalently, $\sigma$ preserves both operations of the associated bi-magmas. The solutions of the braid equation and their morphisms form a category, denoted by $\mathcal{BR}$; in the same way, the YB-solutions form a category $\mathcal{YB}$. As $(\sigma \times \sigma) \circ \tau_X = \tau_{X'} \circ (\sigma \times \sigma)$ for any map $\sigma \colon X \to X'$, a map $\sigma$ is a morphism $(X, \, R) \to (X', \, R')$ if and only if it is a morphism $(X, \, \tau_X \circ R) \to (X', \, \tau_{X'} \circ R')$. Hence the assignment $(X, \, R) \mapsto (X, \, \tau_X \circ R)$, which is the identity on morphisms, defines an isomorphism of categories $\mathcal{BR} \to \mathcal{YB}$, whose inverse is given by the same formula, since $\tau_X^2 = {\rm Id}_{X \times X}$. In particular, if ${\rm Aut} (X, \, R)$ denotes the group of all bijections $\sigma \colon X \to X$ such that $(\sigma \times \sigma) \circ R = R \circ (\sigma \times \sigma)$, in other words the group of automorphisms of the bi-magma $(X, \, \cdot, \, \ast)$, we have
\begin{equation}\eqlabel{auttau}
{\rm Aut} (X, \, R) = {\rm Aut} (X, \, \tau_X \circ R).
\end{equation}
For this reason we shall work only with the braid equation: each result below can be translated into a result on the YB-solutions $r := \tau_X \circ R$. We only have to keep in mind that $R$ is involutive if and only if $r$ is unitary.

Recall that a partially ordered set, or a \emph{poset}, is a set $P$ endowed with a reflexive, antisymmetric and transitive relation $\leq$; we write $x < y$ if $x \leq y$ and $x \neq y$, and we say that $x$ and $y$ are comparable if $x \leq y$ or $y \leq x$. Equivalently, a poset can be given by its strict order relation $<$, which is transitive and irreflexive, i.e. $x \not< x$ for all $x \in P$; this is the point of view of \cite[Chapter IV, 5.1]{PultrTrnkova80} and \cite[Section 2]{Schnyder89}. The two descriptions are equivalent: $\leq$ is obtained from $<$ by adding the equality, and $<$ is obtained from $\leq$ by removing it. A poset $P$ is a \emph{chain} (or a totally ordered set) if any two elements of $P$ are comparable, and an \emph{antichain} if no two distinct elements of $P$ are comparable, i.e. if the order relation $\leq$ is the equality.

The opposite poset $P^{\rm op}$ is the set $P$ with the order given by $x \leq^{\rm op} y$ if and only if $y \leq x$. An automorphism of $P$ is a bijection $\varphi \colon P \to P$ such that $x \leq y$ if and only if $\varphi (x) \leq \varphi (y)$, for all $x$, $y\in P$; equivalently, $x < y$ if and only if $\varphi(x) < \varphi(y)$. The automorphisms of $P$ form a group ${\rm Aut} (P)$. The category of posets and strictly monotone maps, i.e. maps $\sigma$ such that $x < y$ implies $\sigma(x) < \sigma(y)$, is denoted by $\mathcal{P}os$.

Finally, we fix the terminology from graph theory \cite{Bollobas}. A graph $\Gamma = (V, \, E)$ is given by a set $V$, whose elements are the vertices of $\Gamma$, and a set $E$ of subsets of $V$ with two elements, called edges; an edge $\{u, \, v\}$ joins its endpoints $u$ and $v$. The automorphism group ${\rm Aut} (\Gamma)$ is the subgroup of ${\rm Sym} (V)$ consisting of all permutations $\pi$ such that $\{\pi(u), \, \pi(v)\}$ is an edge if and only if $\{u, \, v\}$ is an edge. The graph $\Gamma$ is \emph{bipartite} if $V$ is the disjoint union of two subsets $V_0$ and $V_1$ such that every edge has an endpoint in $V_0$ and the other one
in $V_1$; for a connected bipartite graph, the subsets $V_0$ and $V_1$ are uniquely determined, up to their order (cf. the proof of \cite[Theorem I.4]{Bollobas}). The following theorem was proved by Frucht \cite{Frucht39} for finite groups, and by de Groot \cite{deGroot59} and Sabidussi \cite{Sabidussi60} for arbitrary groups.

\begin{theorem}\thlabel{Frucht}
Any group $G$ is isomorphic to the automorphism group ${\rm Aut} (\Gamma)$ of a graph $\Gamma$, which can be chosen finite if $G$ is finite.
\end{theorem}

%%%%%%%%%%%%%%%%%%%%%%%%%%%%%%%%%%%%%%%%%%%%%%%%%%%%%%%%%%%%%%%%%%%%%%%%%%%%%
%%%%%%%%%%%%%%%%%%%%%%%%%%%%%%%%%%%%%%%%%%%%%%%%%%%%%%%%%%%%%%%%%%%%%%%%%%%%%
\section{Solutions of the braid equation arising from posets}\selabel{posets}

In this section we first extend the example of Lebed \cite[Example 2.7]{Lebed17}, recalled in the introduction, to arbitrary posets, and then we study the new class of solutions obtained in this way from the point of view of Problems A and B of \cite[Introduction]{ACM2}, which are, in our opinion, the basic questions to be asked for any family of solutions of the braid (resp. Yang--Baxter) equation: we prove structure and classification theorems for these solutions, and we characterize them intrinsically, i.e. only in terms of the map $R$.

For a poset $P$ we define $R_P \colon P \times P \to P \times P$ by
\begin{equation}\eqlabel{RP}
R_P (x, \, y) := \begin{cases} (y, \, x), & \text{if} \,\,\, y < x,\\ (x, \, y), & \text{otherwise}. \end{cases}
\end{equation}
Thus, $R_P$ sorts the pair $(x, \, y)$ if $x$ and $y$ are comparable, and leaves it unchanged otherwise; in particular, for all $x$, $y\in P$ we have
\begin{equation}\eqlabel{RPless}
y < x \,\, \Longleftrightarrow \,\, R_P (x, \, y) \neq (x, \, y).
\end{equation}
The associated bi-magma \equref{notatieR} is given by $x \cdot y = y$ and $x \ast y = x$ if $y < x$, and $x \cdot y = x$ and $x \ast y = y$ otherwise. If $P$ is a chain, then $R_P (x, \, y) = (\min\{x, \, y\}, \, \max\{x, \, y\})$ is the solution of the braid equation \cite[Example 2.7]{Lebed17}; if $P$ is an antichain, then $R_P = {\rm Id}_{P \times P}$.

\begin{theorem}\thlabel{RP}
Let $P$ be a poset and $R_P \colon P \times P \to P \times P$ the map given by \equref{RP}. Then:
\begin{enumerate}[(a)]
\item $R_P$ is an idempotent solution of the braid equation.
\item Let $Q$ be another poset. A map $\sigma \colon P \to Q$ is a morphism of solutions $(P, \, R_P) \to (Q, \, R_Q)$ if and only if, for all $x$, $y\in P$:
\begin{equation}\eqlabel{morfRP}
y < x \,\, \Longrightarrow \,\, \sigma (y) \leq \sigma (x) \qquad {\rm and} \qquad \sigma (y) < \sigma (x) \,\, \Longrightarrow \,\, y < x.
\end{equation}
In particular, the isomorphisms of solutions $(P, \, R_P) \to (Q, \, R_Q)$ are exactly the isomorphisms of posets $P \to Q$. Thus, $(P, \, R_P) \cong (Q, \, R_Q)$ in $\BR$ if and only if $P \cong Q$ in the
category $\mathcal{P}os$.
\item ${\rm Aut} (P, \, R_P) = {\rm Aut} (P)$, as subgroups of ${\rm Sym} (P)$.
\end{enumerate}
\end{theorem}

\begin{proof}
(a) Throughout the proof of (a) we write $R := R_P$. If $y < x$, then $R (x, \, y) = (y, \, x)$ and $R (y, \, x) = (y, \, x)$, since $x < y$ is impossible by antisymmetry; otherwise $R (x, \, y) = (x, \, y)$. In both cases $R^2 (x, \, y) = R (x, \, y)$. We show now that $R$ satisfies the braid equation.

Let $(x, \, y, \, z) \in P^3$. We compute the images of $(x, \, y, \, z)$ through the left hand side $R^{12} R^{23} R^{12}$ and through the right hand side $R^{23} R^{12} R^{23}$ of the braid equation, applying the three maps one after another (from right to left). We distinguish four cases, according to whether $y < x$ and whether $z < y$, and we use repeatedly the fact that $R$ acts as the identity on a pair $(u, \, v)$ unless $v < u$.

Case 1: $y \not< x$ and $z \not< y$. Then $R^{12}$ and $R^{23}$ both fix $(x, \, y, \, z)$, hence both sides of the braid equation fix $(x, \, y, \, z)$.

Case 2: $y < x$ and $z \not< y$. The left hand side $R^{12} R^{23} R^{12}$ maps $(x, \, y, \, z)$ as follows: first $(x, \, y, \, z) \overset{R^{12}}{\longmapsto} (y, \, x, \, z)$, and then:
\begin{itemize}
\item if $z < x$: $(y, \, x, \, z) \overset{R^{23}}{\longmapsto} (y, \, z, \, x) \overset{R^{12}}{\longmapsto} (y, \, z, \, x)$, since $z \not< y$;
\item if $z \not< x$: $(y, \, x, \, z) \overset{R^{23}}{\longmapsto} (y, \, x, \, z) \overset{R^{12}}{\longmapsto} (y, \, x, \, z)$, since $x \not< y$.
\end{itemize}
The right hand side $R^{23} R^{12} R^{23}$ maps $(x, \, y, \, z)$ as follows: $(x, \, y, \, z) \overset{R^{23}}{\longmapsto} (x, \, y, \, z) \overset{R^{12}}{\longmapsto} (y, \, x, \, z) \overset{R^{23}}{\longmapsto} (y, \, z, \, x)$ if $z < x$, respectively $(y, \, x, \, z)$ if $z \not< x$. Thus the two sides coincide.

Case 3: $y \not< x$ and $z < y$. The left hand side $R^{12} R^{23} R^{12}$ maps $(x, \, y, \, z)$ as follows: $(x, \, y, \, z) \overset{R^{12}}{\longmapsto} (x, \, y, \, z) \overset{R^{23}}{\longmapsto} (x, \, z, \, y)$, and then $\overset{R^{12}}{\longmapsto} (z, \, x, \, y)$ if $z < x$, respectively $(x, \, z, \, y)$ if $z \not< x$. The right hand side $R^{23} R^{12} R^{23}$ maps $(x, \, y, \, z)$ as follows: first $(x, \, y, \, z) \overset{R^{23}}{\longmapsto} (x, \, z, \, y)$, and then:
\begin{itemize}
\item if $z < x$: $(x, \, z, \, y) \overset{R^{12}}{\longmapsto} (z, \, x, \, y) \overset{R^{23}}{\longmapsto} (z, \, x, \, y)$, since $y \not< x$;
\item if $z \not< x$: $(x, \, z, \, y) \overset{R^{12}}{\longmapsto} (x, \, z, \, y) \overset{R^{23}}{\longmapsto} (x, \, z, \, y)$, since $y \not< z$.
\end{itemize}
Thus the two sides coincide.

Case 4: $y < x$ and $z < y$. By transitivity, $z < x$, so $z < y < x$. The left hand side $R^{12} R^{23} R^{12}$ maps $(x, \, y, \, z)$ as follows:
\begin{equation*}
(x, \, y, \, z) \overset{R^{12}}{\longmapsto} (y, \, x, \, z) \overset{R^{23}}{\longmapsto} (y, \, z, \, x) \overset{R^{12}}{\longmapsto} (z, \, y, \, x),
\end{equation*}
and the right hand side $R^{23} R^{12} R^{23}$ maps it as follows:
\begin{equation*}
(x, \, y, \, z) \overset{R^{23}}{\longmapsto} (x, \, z, \, y) \overset{R^{12}}{\longmapsto} (z, \, x, \, y) \overset{R^{23}}{\longmapsto} (z, \, y, \, x).
\end{equation*}
Thus, in all four cases, the two sides of the braid equation coincide, i.e. $R = R_P$ is a solution of the braid equation. We mention that the transitivity of the order is used only in Case 4.

(b) Let $\sigma \colon P \to Q$ be a map. Then $\sigma$ is a morphism of solutions if and only if
\begin{equation}\eqlabel{morfRPpf}
(\sigma \times \sigma) \bigl( R_P (x, \, y) \bigr) = R_Q \bigl( \sigma(x), \, \sigma(y) \bigr), \qquad \text{for all} \,\, x, \, y\in P.
\end{equation}
We first compute the right hand side: by \equref{RP} applied to $Q$, we have
\begin{equation*}
R_Q \bigl( \sigma(x), \, \sigma(y) \bigr) = \begin{cases} \bigl( \sigma(y), \, \sigma(x) \bigr), & \text{if} \,\,\, \sigma(y) < \sigma(x),\\ \bigl( \sigma(x), \, \sigma(y) \bigr), & \text{otherwise}. \end{cases}
\end{equation*}
For the left hand side of \equref{morfRPpf} we have two cases. In the first case, $y < x$, we have $R_P (x, \, y) = (y, \, x)$ and the left hand side of \equref{morfRPpf} is $(\sigma(y), \, \sigma(x))$. Hence \equref{morfRPpf} holds if and only if $\sigma(y) < \sigma(x)$ or $\sigma(x) = \sigma(y)$, i.e. $\sigma(y) \leq \sigma(x)$.

In the second case, $y \not< x$, we have $R_P (x, \, y) = (x, \, y)$ and the left hand side of \equref{morfRPpf} is $(\sigma(x), \, \sigma(y))$. Hence \equref{morfRPpf} holds if and only if $\sigma(y) \not< \sigma(x)$.

Thus, $\sigma$ is a morphism of solutions if and only if $y < x$ implies $\sigma(y) \leq \sigma(x)$, and $y \not< x$ implies $\sigma(y) \not< \sigma(x)$. The second implication is equivalent to $\sigma(y) < \sigma(x) \Rightarrow y < x$, and we obtain \equref{morfRP}.

We describe now the isomorphisms. Let $\sigma \colon P \to Q$ be a bijection which satisfies \equref{morfRP}. If $y < x$, then $\sigma(y) \leq \sigma(x)$ and $\sigma(y) \neq \sigma(x)$, since $\sigma$ is injective; thus $\sigma(y) < \sigma(x)$. Together with the second implication of \equref{morfRP}, we obtain
\begin{equation*}
y < x \,\, \Longleftrightarrow \,\, \sigma(y) < \sigma(x), \qquad \text{for all} \,\, x, \, y\in P,
\end{equation*}
i.e. $\sigma$ is an isomorphism of posets. Conversely, any isomorphism of posets obviously satisfies \equref{morfRP}. Hence the isomorphisms of solutions $(P, \, R_P) \to (Q, \, R_Q)$ are exactly the isomorphisms of posets $P \to Q$; in particular, $(P, \, R_P) \cong (Q, \, R_Q)$ if and only if $P \cong Q$.

(c) By (b), applied for $Q = P$, the automorphisms of the solution $(P, \, R_P)$ are exactly the automorphisms of the poset $P$, i.e. ${\rm Aut} (P, \, R_P) = {\rm Aut} (P)$.
\end{proof}

The conditions \equref{morfRP} are weaker than the ones defining an order embedding: for instance, any constant map $P \to Q$ is a morphism of solutions $(P, \, R_P) \to (Q, \, R_Q)$. Thus, $P \mapsto (P, \, R_P)$ is the full embedding \equref{PosinYB} of the category $\mathcal{P}os^{\bullet}$ of posets, with the maps satisfying \equref{morfRP} as morphisms, into the category $\BR$ of solutions of the braid equation.

Many families of solutions of the braid equation, or equivalently of the Yang--Baxter equation, are obtained from algebraic or combinatorial structures of various kinds. Once such a family is constructed, a natural question is to recognize its members up to isomorphism: given an arbitrary map $R \colon X \times X \to X \times X$, decide, only by looking at $R$, whether $(X, \, R)$ is isomorphic to a solution of the family. This is Problem B of \cite[Introduction]{ACM2}; in our opinion it is a basic question, although it has been answered for very few families of solutions. For the solutions $R_P$, the answer is given by the following consequence of \thref{RP}:

\begin{corollary}\colabel{caractRP}
Let $X$ be a set and $R \colon X \times X \to X \times X$ a map. Then $(X, \, R)$ is isomorphic to $(P, \, R_P)$, for some poset $P$, if and only if $R$ is an idempotent solution of the braid equation and
\begin{equation}\eqlabel{permut}
R (x, \, y) \in \bigl\{ (x, \, y), \, (y, \, x) \bigr\}, \qquad \text{for all} \,\, x, \, y\in X.
\end{equation}
Moreover, in this case there exists a unique partial order on $X$ such that $R = R_{P_R}$, where $P_R$ denotes the set $X$ with this order, namely
\begin{equation}\eqlabel{orderR}
y < x \,\, \Longleftrightarrow \,\, x \neq y \,\, \text{and} \,\, R (x, \, y) = (y, \, x).
\end{equation}
In particular, $P \mapsto (P, \, R_P)$ induces a bijection between the isomorphism classes of posets with $n$ elements and the isomorphism classes of idempotent solutions of the braid equation on a set with $n$ elements which satisfy \equref{permut}.
\end{corollary}

\begin{proof}
Assume first that $\sigma \colon X \to P$ is an isomorphism of solutions $(X, \, R) \to (P, \, R_P)$, for some poset $P$. Then
\begin{equation*}
R = (\sigma \times \sigma)^{-1} \circ R_P \circ (\sigma \times \sigma).
\end{equation*}
By \thref{RP}(a), $R_P$ is an idempotent solution of the braid equation, hence so is $R$. Moreover, $R_P$ satisfies \equref{permut}, hence $R (x, \, y) = (\sigma \times \sigma)^{-1} \bigl( R_P (\sigma(x), \, \sigma(y)) \bigr) \in \{(x, \, y), \, (y, \, x)\}$, for all $x$, $y\in X$.

Conversely, let $R$ be an idempotent solution of the braid equation which satisfies \equref{permut}, and let $<$ be the relation on $X$ given by \equref{orderR}. By \equref{permut}, for $x \neq y$ we have $R (x, \, y) = (y, \, x)$ if $y < x$, and $R (x, \, y) = (x, \, y)$ otherwise. We show that $<$ is a strict partial order on $X$.

By definition, $x \not< x$, for all $x\in X$. We show that $<$ is asymmetric, i.e. $y < x$ implies $x \not< y$, for all $x$, $y \in X$. If $y < x$, then, as $R$ is idempotent,
\begin{equation*}
R (y, \, x) = R \bigl( R (x, \, y) \bigr) = R (x, \, y) = (y, \, x) \neq (x, \, y),
\end{equation*}
hence $x \not< y$.

We show now that $<$ is transitive. Let $x$, $y$, $z\in X$ be such that $y < x$ and $z < y$; then $z \neq x$, since $x < y$ is impossible by asymmetry. Assume that $z \not< x$. Since $y < x$, $z < y$ and, by asymmetry, $x \not< y$ and $y \not< z$, the two sides of the braid equation give
\begin{align*}
(x, \, y, \, z) &\overset{R^{12}}{\longmapsto} (y, \, x, \, z) \overset{R^{23}}{\longmapsto} (y, \, x, \, z) \overset{R^{12}}{\longmapsto} (y, \, x, \, z),\\
(x, \, y, \, z) &\overset{R^{23}}{\longmapsto} (x, \, z, \, y) \overset{R^{12}}{\longmapsto} (x, \, z, \, y) \overset{R^{23}}{\longmapsto} (x, \, z, \, y).
\end{align*}
These two triples are different, since $x \neq y$, a contradiction. Hence $z < x$, and $<$ is transitive.

Thus, $\leq \, := \, < \cup =$ is a partial order on $X$. If $P_R$ is the set $X$ with this order, then $R = R_{P_R}$, by \equref{RP} and the description of $R$ given above; in particular, $(X, \, R) \cong (P_R, \, R_{P_R})$. If $\leq'$ is another partial order on $X$ such that $R = R_{P'}$, where $P' := (X, \, \leq')$, then by \equref{RPless} we have $y <' x$ if and only if $R (x, \, y) \neq (x, \, y)$, i.e. if and only if $y < x$; hence the order is unique. The last statement follows from the first part and \thref{RP}(b).
\end{proof}

\begin{remarks}\relabel{counting}
\begin{enumerate}[(1)]
\item Let $X$ be a set with $n$ elements. By \coref{caractRP}, the number of idempotent solutions $R$ of the braid equation on $X$ such that $R (x, \, y) \in \{(x, \, y), \, (y, \, x)\}$, for all $x$, $y\in X$, is equal to the number $P_n$ of partial orders on $X$ (i.e. of labelled posets with $n$ elements), and the number of their isomorphism classes is equal to the number $p_n$ of posets with $n$ elements, up to isomorphism. For $n = 1, \ldots, 5$ we have $P_n = 1, \, 3, \, 19, \, 219, \, 4231$ and $p_n = 1, \, 2, \, 5, \, 16, \, 63$ (sequences A001035 and A000112 of \cite{OEIS}). No closed formula is known for $P_n$ or $p_n$.
\item For any poset $P$ we have $\tau_P \circ R_P \circ \tau_P = R_{P^{\rm op}}$. Hence, by \thref{RP}(b), the solutions $R_P$ and $\tau_P \circ R_P \circ \tau_P$ are isomorphic if and only if $P \cong P^{\rm op}$, although they always have the same automorphism group ${\rm Aut} (P) = {\rm Aut} (P^{\rm op})$. For instance, if $P = \{a, \, b, \, c\}$ with $a < b$, $a < c$ and $b$, $c$ incomparable, then $P$ has two maximal elements and $P^{\rm op}$ only one, so $(P, \, R_P) \not\cong (P, \, \tau_P \circ R_P \circ \tau_P)$.
\end{enumerate}
\end{remarks}

By a theorem of Birkhoff \cite[Teorema 2]{Birkhoff46}, every group $G$ is isomorphic to the automorphism group of a poset with $|G|^2 + |G|$ elements. For the convenience of the reader we give a short proof based on \thref{Frucht}. For a graph $\Gamma = (V, \, E)$, its \emph{incidence poset} is the set $P_\Gamma := V \sqcup E$ with the order whose strict relations are $v < e$, for all $v \in V$ and $e \in E$ with $v \in e$ \cite[Section 2]{Schnyder89}; equivalently, $x \leq y$ if and only if either $x = y$, or $x \in V$, $y \in E$ and $x \in y$. This relation is a partial order: it is obviously reflexive and antisymmetric, and it is transitive since there are no elements $a < b < c$ (an element which is strictly greater than another one is an edge, and no element is strictly greater than an edge).

\begin{lemma}\lelabel{incposet}
For any graph $\Gamma = (V, \, E)$, the map ${\rm Aut} (\Gamma) \to {\rm Aut} (P_\Gamma)$, $\pi \mapsto \widehat{\pi}$, where $\widehat{\pi} (v) := \pi (v)$ and $\widehat{\pi} (\{u, \, v\}) := \{\pi(u), \, \pi(v)\}$, is an isomorphism of groups.
\end{lemma}

\begin{proof}
If $\pi \in {\rm Aut} (\Gamma)$, then $\widehat{\pi}$ is a bijection of $P_\Gamma$ mapping $V$ onto $V$ and $E$ onto $E$, and $v \in e$ if and only if $\pi(v) \in \widehat{\pi}(e)$; thus $\widehat{\pi} \in {\rm Aut} (P_\Gamma)$, and $\pi \mapsto \widehat{\pi}$ is clearly an injective morphism of groups. Conversely, let $\varphi \in {\rm Aut} (P_\Gamma)$. The elements of $E$ are exactly the elements of $P_\Gamma$ which are strictly greater than some element (each edge is greater than its two endpoints, and no vertex is greater than any element); this property is preserved by $\varphi$ and by $\varphi^{-1}$, so $\varphi (E) = E$ and $\varphi (V) = V$. Let $\pi := \varphi|_V$. For an edge $e = \{u, \, v\}$, the elements strictly smaller than $e$ are exactly $u$ and $v$; hence the elements strictly smaller than $\varphi(e)$ are exactly $\pi(u)$ and $\pi(v)$, i.e. $\varphi(e) = \{\pi(u), \, \pi(v)\}$. Thus $\pi$ maps edges to edges, and the same argument applied to $\varphi^{-1}$ shows that $\pi^{-1}$ maps edges to edges. Therefore $\pi \in {\rm Aut} (\Gamma)$ and $\varphi = \widehat{\pi}$.
\end{proof}

\begin{corollary}\colabel{poset}
Any group $G$ is isomorphic to the automorphism group of an idempotent solution of the braid equation of the form $R_P$, and to the automorphism group of the YB-solution $\tau_P \circ R_P$, for some poset $P$. If $G \cong {\rm Aut} (\Gamma)$ for a finite graph $\Gamma$ with $n$ vertices and $m$ edges, we can take $P = P_\Gamma$, with $n + m$ elements.
\end{corollary}

\begin{proof}
Let $G$ be a group. By \thref{Frucht}, there exists a graph $\Gamma = (V, \, E)$ such that $G \cong {\rm Aut} (\Gamma)$, and $\Gamma$ can be chosen finite if $G$ is finite. By \leref{incposet}, ${\rm Aut} (\Gamma) \cong {\rm Aut} (P_\Gamma)$, where $P_\Gamma$ is the incidence poset of $\Gamma$, which has $|V| + |E|$ elements. Finally, by \thref{RP}(a) and (c), $R_{P_\Gamma}$ is an idempotent solution of the braid equation and ${\rm Aut} (P_\Gamma, \, R_{P_\Gamma}) = {\rm Aut} (P_\Gamma)$; moreover, ${\rm Aut} (P_\Gamma, \, \tau_{P_\Gamma} \circ R_{P_\Gamma}) = {\rm Aut} (P_\Gamma, \, R_{P_\Gamma})$, by \equref{auttau}. Therefore,
\begin{equation*}
G \cong {\rm Aut} (\Gamma) \cong {\rm Aut} (P_\Gamma) = {\rm Aut} (P_\Gamma, \, R_{P_\Gamma}) = {\rm Aut} (P_\Gamma, \, \tau_{P_\Gamma} \circ R_{P_\Gamma}),
\end{equation*}
and the poset $P := P_\Gamma$ has the required properties.
\end{proof}

\begin{remark}\relabel{poset}
Let $L$ be a lattice and $P$ its underlying poset. Then $R_P (x, \, y) = (x \wedge y, \, x \vee y)$ if $x$ and $y$ are comparable, and $R_P (x, \, y) = (x, \, y)$ otherwise. On the other hand, the map $R_L \colon L \times L \to L \times L$, $R_L (x, \, y) := (x \wedge y, \, x \vee y)$, is a solution of the braid equation if and only if $L$ is distributive: see \cite[Example 2.8]{Lebed17} and \cite[Section 4 and Proposition 5.13]{CVV20}. The solutions $R_L$ and $R_P$ coincide if and only if $L$ is a chain, and in this case they are equal to the solution of \cite[Example 2.7]{Lebed17}.
\end{remark}

The idempotent solutions $R$ of the braid equation such that $R(x, \, y) \in \{(x, \, y), \, (y, \, x)\}$, for all $x$, $y$, are described by \coref{caractRP}. We drop now the idempotency assumption.
Let $X$ be a set. Any map $R \colon X \times X \to X \times X$ with this property is determined by the set of pairs which it interchanges,
\begin{equation}\eqlabel{SR}
S_R := \bigl\{ (x, \, y) \in X \times X \, \bigl| \, x \neq y, \,\, R (x, \, y) = (y, \, x) \bigr\},
\end{equation}
since $R (x, \, x) = (x, \, x)$, for all $x \in X$. Conversely, let $S \subseteq X \times X$ be a subset disjoint from the diagonal, i.e. $(x, \, x) \notin S$, for all $x\in X$. Then the map $R \colon X \times X \to X \times X$ defined by $R (x, \, y) := (y, \, x)$ if $(x, \, y) \in S$, and $R (x, \, y) := (x, \, y)$ otherwise, satisfies $R (x, \, y) \in \{(x, \, y), \, (y, \, x)\}$ and $S_R = S$, and it is the only map with these properties. Thus, $R \mapsto S_R$ is a bijection between the maps $R$ with $R(x, \, y) \in \{(x, \, y), \, (y, \, x)\}$, for all $x$, $y\in X$, and the subsets of $X \times X$ disjoint from the diagonal.

As a first step towards the description of this family, we generalize the solutions $R_P$ given by \equref{RP} as follows:

\begin{definition}\delabel{partinpos}
A \emph{decomposition of a set $X$ into posets} is a family $\mathcal{P} = (X_i, \, \leq_i)_{i\in I}$, where $(X_i)_{i\in I}$ is a partition of $X = \sqcup_{i\in I} \, X_i$ into non-empty subsets, called the blocks of $\mathcal{P}$, and $\leq_i$ is a partial order on $X_i$, for all $i\in I$. To such a decomposition we associate the map
\begin{equation}\eqlabel{RcalP}
R_{\mathcal{P}} \colon X \times X \to X \times X, \quad R_{\mathcal{P}} (x, \, y) := \begin{cases} R_{X_i} (x, \, y), & \text{if} \,\,\, x, \, y \in X_i, \,\, \text{for some} \,\, i\in I,\\ (y, \, x), & \text{if} \,\,\, x \,\, \text{and} \,\, y \,\, \text{belong to different blocks}, \end{cases}
\end{equation}
where $R_{X_i}$ is the map \equref{RP} of the poset $(X_i, \, \leq_i)$.
\end{definition}

Thus, $R_{\mathcal{P}}$ sorts the comparable pairs inside each block, fixes the incomparable ones, and flips the pairs whose components belong to different blocks. If $\mathcal{P}$ has only one
block we recover $R_P$; if all the blocks are antichains we obtain the map
\begin{equation}\eqlabel{Rpi}
R_\pi (x, \, y) := \begin{cases} (y, \, x), & \text{if} \,\,\, x \,\, \text{and} \,\, y \,\, \text{belong to different blocks},\\ (x, \, y), & \text{otherwise}, \end{cases}
\end{equation}
which depends only on the partition $\pi := (X_i)_{i\in I}$ of $X$. For a partition with two blocks $Y$ and $Z$, $R_\pi$ is the bipartite solution associated to the split datum $(Y, \, Z, \, ({\rm Id}_Z)_{y\in Y})$
(see \cite[Definition 2.1]{MilBIP}).

\begin{theorem}\thlabel{flipfix}
Let $X$ be a set and $R \colon X \times X \to X \times X$ a map such that $R (x, \, y) \in \{(x, \, y), \, (y, \, x)\}$, for all $x$, $y\in X$. Then:
\begin{enumerate}[(a)]
\item $R$ is a solution of the braid equation if and only if $R = R_{\mathcal{P}}$, for some decomposition $\mathcal{P} = (X_i, \, \leq_i)_{i\in I}$ of $X$ into posets. In this case $\mathcal{P}$ is unique: two distinct elements $x$, $y \in X$ belong to different blocks if and only if $(x, \, y) \in S_R$ and $(y, \, x) \in S_R$, while for $x$, $y$ in the same block $X_i$ we have $y <_i x$ if and only if $(x, \, y) \in S_R$;
\item $R_{\mathcal{P}}^2 = R_{\mathcal{P}}$ if and only if $\mathcal{P}$ has at most one block;
\item $R_{\mathcal{P}}$ is bijective if and only if all the blocks of $\mathcal{P}$ are antichains, i.e. $R_{\mathcal{P}} = R_\pi$;
\item if $\mathcal{Q} = (X'_j, \, \leq'_j)_{j\in J}$ is a decomposition of a set $X'$ into posets, then a bijection $\sigma \colon X \to X'$ is an isomorphism of solutions $(X, \, R_{\mathcal{P}}) \to (X', \, R_{\mathcal{Q}})$ if and only if there exists a bijection $\theta \colon I \to J$ such that $\sigma (X_i) = X'_{\theta(i)}$ and $\sigma|_{X_i} \colon (X_i, \, \leq_i) \to (X'_{\theta(i)}, \, \leq'_{\theta(i)})$ is an isomorphism of posets, for all $i\in I$. In particular, ${\rm Aut} (X, \, R_{\mathcal{P}})$ is the group of all bijections of $X$ which permute the blocks of $\mathcal{P}$ and restrict to isomorphisms of posets between the corresponding blocks.
\end{enumerate}
\end{theorem}

\begin{proof}
(a) Assume first that $R = R_{\mathcal{P}}$ and let $(x, \, y, \, z) \in X^3$. If $x$, $y$, $z$ belong to the same block, the braid equation holds by \thref{RP}(a). If they belong to three different blocks, then $R^{12}$ and $R^{23}$ act as the flip on all the pairs which occur, and both sides of the braid equation map $(x, \, y, \, z)$ to $(z, \, y, \, x)$. In the remaining three cases, let $X_i \neq X_j$ be two blocks, and for $u$, $v \in X_i$ write $R_{X_i} (u, \, v) = (u', \, v')$. If $u$, $v\in X_i$ and $w \in X_j$, the two sides of the braid equation give:
\begin{align*}
(u, \, v, \, w) &\overset{R^{12}}{\longmapsto} (u', \, v', \, w) \overset{R^{23}}{\longmapsto} (u', \, w, \, v') \overset{R^{12}}{\longmapsto} (w, \, u', \, v'),\\
(u, \, v, \, w) &\overset{R^{23}}{\longmapsto} (u, \, w, \, v) \overset{R^{12}}{\longmapsto} (w, \, u, \, v) \overset{R^{23}}{\longmapsto} (w, \, u', \, v');\\[4pt]
(w, \, u, \, v) &\overset{R^{12}}{\longmapsto} (u, \, w, \, v) \overset{R^{23}}{\longmapsto} (u, \, v, \, w) \overset{R^{12}}{\longmapsto} (u', \, v', \, w),\\
(w, \, u, \, v) &\overset{R^{23}}{\longmapsto} (w, \, u', \, v') \overset{R^{12}}{\longmapsto} (u', \, w, \, v') \overset{R^{23}}{\longmapsto} (u', \, v', \, w);\\[4pt]
(u, \, w, \, v) &\overset{R^{12}}{\longmapsto} (w, \, u, \, v) \overset{R^{23}}{\longmapsto} (w, \, u', \, v') \overset{R^{12}}{\longmapsto} (u', \, w, \, v'),\\
(u, \, w, \, v) &\overset{R^{23}}{\longmapsto} (u, \, v, \, w) \overset{R^{12}}{\longmapsto} (u', \, v', \, w) \overset{R^{23}}{\longmapsto} (u', \, w, \, v').
\end{align*}
Thus $R_{\mathcal{P}}$ is a solution of the braid equation.

Conversely, let $R$ be a solution of the braid equation and $S := S_R$. We write $x \rightleftarrows y$ if $(x, \, y) \in S$ and $(y, \, x) \in S$, and $x \approx y$ if $x \rightleftarrows y$ does not hold; the relation $\approx$ is reflexive, since $S$ is disjoint from the diagonal, and symmetric. We show that $\approx$ is transitive. Assume, by contradiction, that $x \approx z$, $z \approx y$ and $x \rightleftarrows y$; then $x \neq y$ and $z \notin \{x, \, y\}$. We distinguish three cases.

Case 1: $(x, \, z) \notin S$. The left hand side of the braid equation gives $(x, \, y, \, z) \mapsto (y, \, x, \, z) \mapsto (y, \, x, \, z) \mapsto (x, \, y, \, z)$. If $(y, \, z) \in S$, then $(z, \, y) \notin S$, since $z \approx y$, and the right hand side gives $(x, \, y, \, z) \mapsto (x, \, z, \, y) \mapsto (x, \, z, \, y) \mapsto (x, \, z, \, y) \neq (x, \, y, \, z)$; if $(y, \, z) \notin S$, the right hand side gives $(x, \, y, \, z) \mapsto (x, \, y, \, z) \mapsto (y, \, x, \, z) \mapsto (y, \, x, \, z) \neq (x, \, y, \, z)$. In both cases we obtain a contradiction.

Case 2: $(y, \, z) \notin S$. This case follows from Case 1 by interchanging $x$ and $y$.

Case 3: $(x, \, z) \in S$ and $(y, \, z) \in S$. Then $(z, \, x) \notin S$ and $(z, \, y) \notin S$, since $z \approx x$ and $z \approx y$. The left hand side gives $(z, \, x, \, y) \mapsto (z, \, x, \, y) \mapsto (z, \, y, \, x) \mapsto (z, \, y, \, x)$, while the right hand side gives $(z, \, x, \, y) \mapsto (z, \, y, \, x) \mapsto (z, \, y, \, x) \mapsto (z, \, x, \, y)$, a contradiction.

Thus $\approx$ is an equivalence relation on $X$; let $(X_i)_{i\in I}$ be its classes. If $x$ and $y$ belong to different classes, then $x \rightleftarrows y$, i.e. $R (x, \, y) = (y, \, x)$. For $x$, $y$ in the same class $X_i$ we define $y <_i x$ if and only if $(x, \, y) \in S$. The relation $<_i$ is irreflexive, and it is asymmetric since $X_i$ contains no pair $x \rightleftarrows y$. It is also transitive: let $z <_i y <_i x$, so that $(x, \, y)$, $(y, \, z) \in S$ and $(y, \, x)$, $(z, \, y) \notin S$; by asymmetry, $z \neq x$. If $(x, \, z) \notin S$, then the left hand side of the braid equation gives $(x, \, y, \, z) \mapsto (y, \, x, \, z) \mapsto (y, \, x, \, z) \mapsto (y, \, x, \, z)$, while the right hand side gives $(x, \, y, \, z) \mapsto (x, \, z, \, y) \mapsto (x, \, z, \, y) \mapsto (x, \, z, \, y)$, a contradiction. Hence $z <_i x$. Therefore $\leq_i \, := \, <_i \cup =$ is a partial order on $X_i$ and, by construction, $R = R_{\mathcal{P}}$, where $\mathcal{P} := (X_i, \, \leq_i)_{i\in I}$. The uniqueness of $\mathcal{P}$ is clear, since the blocks and the orders of any decomposition $\mathcal{P}'$ with $R = R_{\mathcal{P}'}$ are recovered from $S_R$ as above.

(b) If $x$ and $y$ belong to different blocks, then $R_{\mathcal{P}}^2 (x, \, y) = (x, \, y) \neq (y, \, x) = R_{\mathcal{P}} (x, \, y)$. Hence, if $R_{\mathcal{P}}^2 = R_{\mathcal{P}}$, then $\mathcal{P}$ has at most one block; the converse follows from \thref{RP}(a).

(c) If $y <_i x$ in some block $X_i$, then $R_{\mathcal{P}} (x, \, y) = (y, \, x) = R_{\mathcal{P}} (y, \, x)$, and $R_{\mathcal{P}}$ is not injective. If all the blocks are antichains, then $R_{\mathcal{P}} = R_\pi$ and obviously $R_\pi^2 = {\rm Id}_{X \times X}$.

(d) Let $R$ and $R'$ be maps as in the statement on $X$, respectively $X'$, and $\sigma \colon X \to X'$ a bijection. We show first that $(\sigma \times \sigma) \circ R = R' \circ (\sigma \times \sigma)$ if and only if $(\sigma \times \sigma) (S_R) = S_{R'}$. Let $x$, $y\in X$. If $(x, \, y) \in S_R$, then $\sigma(x) \neq \sigma(y)$ and $(\sigma \times \sigma) (R (x, \, y)) = (\sigma(y), \, \sigma(x))$, which is equal to $R' (\sigma(x), \, \sigma(y))$ if and only if $(\sigma(x), \, \sigma(y)) \in S_{R'}$. If $(x, \, y) \notin S_R$, then $(\sigma \times \sigma) (R (x, \, y)) = (\sigma(x), \, \sigma(y))$, which is equal to $R' (\sigma(x), \, \sigma(y))$ if and only if $(\sigma(x), \, \sigma(y)) \notin S_{R'}$. Now let $R = R_{\mathcal{P}}$ and $R' = R_{\mathcal{Q}}$. The set $S_{R_{\mathcal{P}}}$ is the disjoint union of the set of pairs of elements belonging to different blocks of $\mathcal{P}$, which is invariant under $(x, \, y) \mapsto (y, \, x)$, and of the set of pairs $(x, \, y)$ with $y <_i x$ for some $i$, which contains no pair together with its transpose; the same holds for $S_{R_{\mathcal{Q}}}$. Since $\sigma \times \sigma$ commutes with $(x, \, y) \mapsto (y, \, x)$, it follows that $(\sigma \times \sigma) (S_{R_\mathcal{P}}) = S_{R_\mathcal{Q}}$ if and only if $\sigma$ maps the pairs of elements in different blocks of $\mathcal{P}$ onto the pairs of elements in different blocks of $\mathcal{Q}$, i.e. $\sigma$ maps the blocks of $\mathcal{P}$ onto the blocks of $\mathcal{Q}$, and, for $x$, $y$ in the same block $X_i$, $y <_i x$ if and only if $\sigma(y) <'_{\theta(i)} \sigma(x)$. This proves (d).
\end{proof}

\begin{remark}\relabel{flip}
Let $X$ be a set with $n$ elements. By \thref{flipfix}, the number of solutions $R$ of the braid equation on $X$ such that $R (x, \, y) \in \{(x, \, y), \, (y, \, x)\}$, for all $x$, $y\in X$, is the number of decompositions of $X$ into posets; its exponential generating function is $\exp \bigl( \sum_{k \geq 1} P_k \, t^k / k! \bigr)$, where $P_k$ is the number of labelled posets with $k$ elements (\reref{counting}). For $n = 1, \ldots, 5$ we obtain $1$, $4$, $29$, $341$, $6252$ solutions. By \thref{flipfix}(d), the isomorphism classes of these solutions are in bijection with the multisets (i.e. collections in which repetitions are allowed) of isomorphism classes of non-empty posets whose cardinalities add up to $n$; their number is $1$, $3$, $8$, $27$, $100$, for $n = 1, \ldots, 5$. These numbers were obtained from the values of $P_k$ and $p_k$ given in \cite{OEIS}.
\end{remark}

%%%%%%%%%%%%%%%%%%%%%%%%%%%%%%%%%%%%%%%%%%%%%%%%%%%%%%%%%%%%%%%%%%%%%%%%%%%%%

\section{Two families of involutive solutions associated to posets}\selabel{invol}

By \thref{flipfix}(c), a bijective solution of the braid equation on a poset $P$ which permutes the components of each pair is of the form $R_\pi$, hence it depends only on a partition of $P$. In order to obtain involutive (in particular, bijective) solutions which take into account the order of $P$, we have to enlarge the underlying set; for this purpose we use the bipartite solutions of \cite{MilBIP}. Following \cite[Definition 2.1]{MilBIP}, a \emph{bipartite datum} is a triple $(Y, \, Z, \, g)$, where $Y$ and $Z$ are two disjoint sets and $g = (g_y)_{y\in Y}$ is a family of permutations of $Z$, without any relation between them; the datum is called \emph{split} if $g_y = {\rm Id}_Z$ for all $y\in Y$, and \emph{non-split} otherwise. The associated \emph{bipartite solution} is the map $R_g \colon X \times X \to X \times X$, where $X := Y \sqcup Z$, given by
\begin{equation}\eqlabel{Rg}
R_g (x, \, x') := \begin{cases} (x, \, x'), & \text{if} \,\,\, x, \, x' \in Y \,\,\, \text{or} \,\,\, x, \, x'\in Z,\\
\bigl( g_y (z), \, y \bigr), & \text{if} \,\,\, x = y \in Y, \,\, x' = z \in Z,\\
\bigl( y, \, g_y^{-1} (z) \bigr), & \text{if} \,\,\, x = z \in Z, \,\, x' = y \in Y.
\end{cases}
\end{equation}
By \cite[Theorem 2.2]{MilBIP}, $R_g$ is a solution of the braid equation and $R_g^2 = {\rm Id}_{X \times X}$. Moreover, if the datum is non-split, then by \cite[Proposition 2.15]{MilBIP} the map $H_g \to {\rm Aut} (X, \, R_g)$, $(\alpha, \, \beta) \mapsto \alpha \sqcup \beta$, is an isomorphism of groups, where
\begin{equation}\eqlabel{Hg}
H_g := \Bigl\{ (\alpha, \, \beta) \in {\rm Sym} (Y) \times {\rm Sym} (Z) \, \bigl| \, \beta \circ g_y \circ \beta^{-1} = g_{\alpha(y)}, \,\, \text{for all} \,\, y\in Y \Bigr\}
\end{equation}
and $\alpha \sqcup \beta$ is the bijection of $X$ which is equal to $\alpha$ on $Y$ and to $\beta$ on $Z$.

We shall apply this construction to two bipartite data associated to a poset $P$. For $x\in P$ we denote ${\downarrow} x := \{ y\in P \, | \, y \leq x \}$. The following elementary lemma will be used twice.

\begin{lemma}\lelabel{downsets}
Let $P$ and $Q$ be two finite posets and $\alpha$, $\beta \colon P \to Q$ two bijections such that $\alpha ({\downarrow} x) = {\downarrow} \beta(x)$, for all $x\in P$. Then $\alpha = \beta$ and $\alpha$ is an isomorphism of posets.
\end{lemma}

\begin{proof}
We show that $\alpha(x) = \beta(x)$, for all $x \in P$, by induction on $|{\downarrow} x|$. Let $x\in P$ and assume that $\alpha(w) = \beta(w)$ for all $w < x$; this is the induction hypothesis, since $w < x$ implies ${\downarrow} w \subsetneq {\downarrow} x$, hence $|{\downarrow} w| < |{\downarrow} x|$. Since $\beta(x) \in {\downarrow} \beta(x) = \alpha ({\downarrow} x)$, there exists $w \in {\downarrow} x$, i.e. $w \leq x$, such that $\beta(x) = \alpha(w)$. Assume that $w < x$. Then $\alpha(w) = \beta(w)$, by the induction hypothesis, and hence $\beta(x) = \alpha(w) = \beta(w)$.
As $\beta$ is injective, we obtain $w = x$, which contradicts $w < x$. Thus $w = x$, i.e. $\beta(x) = \alpha(x)$.

Finally, since $\alpha$ is bijective and $\alpha = \beta$, for all $x$, $y \in P$ we have: $y \leq x$ if and only if $y \in {\downarrow} x$, if and only if $\alpha(y) \in \alpha ({\downarrow} x) = {\downarrow} \alpha(x)$, i.e. $\alpha(y) \leq \alpha(x)$. Hence $\alpha$ is an isomorphism of posets.
\end{proof}

Let $P$ be a poset and $Z_P := P \times \bZ/2$, where $\bZ/2 = \{0, \, 1\}$ is the usual additive group. For $x\in P$ let $\chi_x \colon P \to \bZ/2$ be the characteristic function of ${\downarrow} x$, and let $d_x \in {\rm Sym} (Z_P)$ be the involution defined for any $(y, \, \varepsilon) \in Z_P$ by
\begin{equation*}
 d_x \colon Z_P \to Z_P, \qquad   d_x (y, \, \varepsilon) := \bigl( y, \, \varepsilon + \chi_x (y) \bigr),
\end{equation*}
i.e. $d_x$ changes the second component of the elements $(y, \, \varepsilon)$ with $y \leq x$, and fixes the other elements of $Z_P$. We denote by $D_P := R_d$ the bipartite solution associated to the datum $(P, \, Z_P, \, d)$, $d = (d_x)_{x\in P}$, on the set $X_P := P \sqcup Z_P$, which has $3 |P|$ elements if $P$ is finite. Explicitly, by \equref{Rg}, for $x\in P$ and $(y, \, \varepsilon) \in Z_P$ we have
\begin{equation}\eqlabel{DP}
D_P \bigl( x, \, (y, \, \varepsilon) \bigr) = \bigl( (y, \, \varepsilon + \chi_x (y)), \, x \bigr), \qquad D_P \bigl( (y, \, \varepsilon), \, x \bigr) = \bigl( x, \, (y, \, \varepsilon + \chi_x (y)) \bigr),
\end{equation}
while $D_P$ is the identity on $P \times P$ and on $Z_P \times Z_P$. By \cite[Theorem 2.2]{MilBIP}, $D_P \colon X_P \times X_P \to X_P \times X_P$ is a solution of the braid equation and $D_P^2 = {\rm Id}$.

Let $(\bZ/2)^P$ be the additive group of all maps $c \colon P \to \bZ/2$, with the pointwise addition. The group ${\rm Aut} (P)$ acts on $(\bZ/2)^P$ by group automorphisms, via $(\varphi \triangleright c) (u) := c \bigl(\varphi^{-1} (u) \bigr)$, for all $\varphi \in {\rm Aut} (P)$, $c\in (\bZ/2)^P$ and $u\in P$, and we denote by $(\bZ/2)^P \rtimes {\rm Aut} (P)$ the corresponding semidirect product: it is the set $(\bZ/2)^P \times {\rm Aut} (P)$ with the multiplication $(c, \, \varphi) (c', \, \varphi') := (c + \varphi \triangleright c', \, \varphi \varphi')$.

For $c\in (\bZ/2)^P$ and $\varphi \in {\rm Aut} (P)$ we define the bijection $\sigma_{c, \varphi} \colon X_P \to X_P$ by
\begin{equation}\eqlabel{sigmacphi}
\sigma_{c, \varphi} (x) := \varphi (x), \qquad \sigma_{c, \varphi} (y, \, \varepsilon) := \Bigl( \varphi (y), \, \varepsilon + c \bigl( \varphi (y) \bigr) \Bigr),
\end{equation}
for all $x\in P$ and $(y, \, \varepsilon) \in Z_P$.

\begin{theorem}\thlabel{autDP}
Let $P$ be a finite non-empty poset and $D_P$ the involutive solution of the braid equation on $X_P = P \sqcup Z_P = P \sqcup (P\times \bZ/2)$ given by \equref{DP}. Then the map
\begin{equation*}
(\bZ/2)^P \rtimes {\rm Aut} (P) \to {\rm Aut} (X_P, \, D_P), \qquad (c, \, \varphi) \mapsto \sigma_{c, \varphi},
\end{equation*}
is an isomorphism of groups, with $\sigma_{c, \varphi}$ given by \equref{sigmacphi}. In particular, $|{\rm Aut} (X_P, \, D_P)| = 2^{|P|} \, |{\rm Aut} (P)|$, and ${\rm Aut} (P)$ is isomorphic to the quotient of ${\rm Aut} (X_P, \, D_P)$ by the normal subgroup of the automorphisms which fix $P$ pointwise.
\end{theorem}

\begin{proof}
The datum $(P, \, Z_P, \, d)$ is non-split: since $P$ is non-empty, we can pick $x\in P$, and then $d_x (x, \, 0) = (x, \, 1)$, as $x \leq x$; thus $d_x \neq {\rm Id}_{Z_P}$. Hence, by \cite[Proposition 2.15]{MilBIP}, the map $H_d \to {\rm Aut} (X_P, \, D_P)$, $(\alpha, \, \beta) \mapsto \alpha \sqcup \beta$, is an isomorphism of groups, where, by \equref{Hg},
\begin{equation*}
H_d = \Bigl\{ (\alpha, \, \beta) \in {\rm Sym} (P) \times {\rm Sym} (Z_P) \, \bigl| \, \beta \circ d_x \circ \beta^{-1} = d_{\alpha(x)}, \,\, \text{for all} \,\, x\in P \Bigr\}.
\end{equation*}
For $c \in (\bZ/2)^P$ and $\varphi \in {\rm Aut} (P)$ we denote by $\beta_{c, \varphi}$ the restriction of $\sigma_{c, \varphi}$ to $Z_P$, i.e. $\beta_{c, \varphi} (y, \, \varepsilon) = \bigl( \varphi (y), \, \varepsilon + c (\varphi(y)) \bigr)$, so that $\sigma_{c, \varphi} = \varphi \sqcup \beta_{c, \varphi}$. Therefore, in order to describe ${\rm Aut} (X_P, \, D_P)$, it is enough to show that the elements of $H_d$ are exactly the pairs $(\varphi, \, \beta_{c, \varphi})$, with $c\in (\bZ/2)^P$ and $\varphi \in {\rm Aut} (P)$.

We first show that $(\varphi, \, \beta_{c, \varphi}) \in H_d$. For $y \in P$ let $e_y := \{ (y, \, 0), \, (y, \, 1) \}$; then $Z_P$ is the disjoint union of the sets $e_y$, $y \in P$, and $\beta_{c, \varphi}$ maps $e_y$ bijectively onto $e_{\varphi(y)}$, for all $y\in P$. As $\varphi$ is bijective, it follows that $\beta_{c, \varphi} \in {\rm Sym} (Z_P)$. Let now $x\in P$ and $(y, \, \varepsilon) \in Z_P$. On the one hand,
\begin{equation*}
\beta_{c, \varphi} \bigl( d_x (y, \, \varepsilon) \bigr) = \beta_{c, \varphi} \bigl( y, \, \varepsilon + \chi_x (y) \bigr) = \Bigl( \varphi(y), \, \varepsilon + \chi_x (y) + c \bigl( \varphi (y) \bigr) \Bigr),
\end{equation*}
and, on the other hand,
\begin{equation*}
d_{\varphi(x)} \bigl( \beta_{c, \varphi} (y, \, \varepsilon) \bigr) = d_{\varphi(x)} \Bigl( \varphi (y), \, \varepsilon + c \bigl( \varphi(y) \bigr) \Bigr) = \Bigl( \varphi(y), \, \varepsilon + c \bigl( \varphi (y) \bigr) + \chi_{\varphi(x)} \bigl( \varphi(y) \bigr) \Bigr).
\end{equation*}
Since $\varphi$ is an automorphism of $P$, we have $y \leq x$ if and only if $\varphi(y) \leq \varphi(x)$, i.e. $\chi_x (y) = \chi_{\varphi(x)} (\varphi(y))$. Thus $\beta_{c, \varphi} \circ d_x = d_{\varphi(x)} \circ \beta_{c, \varphi}$, for all $x\in P$, i.e. $(\varphi, \, \beta_{c, \varphi}) \in H_d$.

Conversely, let $(\alpha, \, \beta) \in H_d$. We show that $\alpha$ is an automorphism of $P$ and $\beta = \beta_{c, \alpha}$, for some $c \in (\bZ/2)^P$. By definition, $d_x$ interchanges the two elements of $e_y$ if $y \leq x$, and fixes them if $y \not\leq x$. Thus the elements of $Z_P$ moved by $d_x$ are exactly the elements of $\bigcup_{y \leq x} e_y$, and $d_x$ interchanges the two elements of each $e_y$ with $y \leq x$. On the other hand, the equality $\beta \circ d_x \circ \beta^{-1} = d_{\alpha(x)}$ means that
\begin{equation*}
\beta \bigl( d_x (z) \bigr) = d_{\alpha(x)} \bigl( \beta (z) \bigr), \qquad \text{for all} \,\, x \in P \,\, \text{and} \,\, z\in Z_P.
\end{equation*}
Hence $\beta$ maps the set of elements moved by $d_x$ onto the set of elements moved by $d_{\alpha(x)}$, i.e.
\begin{equation}\eqlabel{betadown}
\beta \Bigl( \bigcup_{y \leq x} e_y \Bigr) = \bigcup_{w \leq \alpha(x)} e_w, \qquad \text{for all} \,\, x\in P,
\end{equation}
and it maps any pair $\{z, \, d_x (z)\}$ with $z \neq d_x (z)$ onto the pair $\{\beta(z), \, d_{\alpha(x)} (\beta(z))\}$. Applying the last property for $x = y$ and $z = (y, \, 0)$, we obtain that $\beta (e_y)$ is a pair of elements interchanged by $d_{\alpha(y)}$, i.e. $\beta (e_y) = e_{\bar{\beta}(y)}$ for some $\bar{\beta}(y) \in P$. The map $\bar{\beta} \colon P \to P$ obtained in this way is injective, since the sets $e_y$ are pairwise disjoint and $\beta$ is injective; as $P$ is finite, $\bar{\beta}$ is bijective. Now \equref{betadown} can be written as $\bar{\beta} ({\downarrow} x) = {\downarrow} \alpha (x)$, for all $x\in P$. By \leref{downsets}, applied to the bijections $\bar{\beta}$ and $\alpha$, we obtain that $\bar{\beta} = \alpha$ and that $\varphi := \alpha$ is an automorphism of $P$.

Thus, for any $y\in P$, $\beta$ restricts to a bijection $e_y \to e_{\varphi (y)}$. There are only two such bijections, namely $(y, \, \varepsilon) \mapsto (\varphi(y), \, \varepsilon)$ and $(y, \, \varepsilon) \mapsto (\varphi(y), \, \varepsilon + 1)$; hence there exists a unique $c_y \in \bZ/2$ such that $\beta (y, \, \varepsilon) = (\varphi (y), \, \varepsilon + c_y)$, for $\varepsilon \in \bZ/2$. Let $c \in (\bZ/2)^P$ be given by $c (\varphi(y)) := c_y$, for all $y\in P$; this is well defined, since $\varphi$ is bijective. Then $\beta = \beta_{c, \varphi}$, and $\alpha \sqcup \beta = \sigma_{c, \varphi}$.

We have proved that ${\rm Aut} (X_P, \, D_P) = \bigl\{ \sigma_{c, \varphi} \, | \, c\in (\bZ/2)^P, \,\, \varphi \in {\rm Aut} (P) \bigr\}$, i.e. the map $(c, \, \varphi) \mapsto \sigma_{c, \varphi}$ is surjective. It is also injective: if $\sigma_{c, \varphi} = \sigma_{c', \varphi'}$, then $\varphi = \varphi'$, by restricting to $P$, and then, comparing the second components of $\sigma_{c, \varphi} (y, \, 0)$ and $\sigma_{c', \varphi} (y, \, 0)$, we obtain $c (\varphi (y)) = c' (\varphi(y))$, for all $y\in P$, i.e. $c = c'$. Finally, it is a morphism of groups. Indeed, let $(c, \, \varphi)$, $(c', \, \varphi') \in (\bZ/2)^P \rtimes {\rm Aut} (P)$. Both $\sigma_{c, \varphi} \circ \sigma_{c', \varphi'}$ and $\sigma_{c + \varphi \triangleright c', \, \varphi \varphi'}$ are equal to $\varphi \varphi'$ on $P$. For $(y, \, \varepsilon) \in Z_P$ we have
\begin{equation*}
\sigma_{c, \varphi} \bigl( \sigma_{c', \varphi'} (y, \, \varepsilon) \bigr) = \sigma_{c, \varphi} \Bigl( \varphi' (y), \, \varepsilon + c' \bigl( \varphi'(y) \bigr) \Bigr) = \Bigl( \varphi \varphi'(y), \, \varepsilon + c' \bigl( \varphi'(y) \bigr) + c \bigl( \varphi \varphi' (y) \bigr) \Bigr),
\end{equation*}
while, since $(\varphi \triangleright c') (\varphi \varphi'(y)) = c' (\varphi' (y))$,
\begin{equation*}
\sigma_{c + \varphi \triangleright c', \, \varphi \varphi'} (y, \, \varepsilon) = \Bigl( \varphi \varphi'(y), \, \varepsilon + c \bigl( \varphi \varphi'(y) \bigr) + c' \bigl( \varphi'(y) \bigr) \Bigr).
\end{equation*}
Thus $\sigma_{c, \varphi} \circ \sigma_{c', \varphi'} = \sigma_{(c, \varphi) (c', \varphi')}$, and the first statement of the theorem is proved.

Finally, the map $(\bZ/2)^P \rtimes {\rm Aut} (P) \to {\rm Aut} (P)$, $(c, \, \varphi) \mapsto \varphi$, is a surjective morphism of groups, whose kernel is $(\bZ/2)^P \times \{ {\rm Id}_P \}$. Under the isomorphism proved above, this kernel corresponds to the set of automorphisms $\sigma_{c, {\rm Id}_P}$, $c \in (\bZ/2)^P$, which are exactly the automorphisms of $(X_P, \, D_P)$ fixing $P$ pointwise, since $\sigma_{c, \varphi}$ is equal to $\varphi$ on $P$. Hence these automorphisms form a normal subgroup of ${\rm Aut} (X_P, \, D_P)$, and the quotient by it is isomorphic to ${\rm Aut} (P)$.
\end{proof}

\begin{corollary}\colabel{isoDP}
Let $P$ and $Q$ be finite non-empty posets. Then $(X_P, \, D_P) \cong (X_Q, \, D_Q)$ if and only if $P \cong Q$.\footnote{This answers a question addressed to the author in a private 
communication by Arhan Canli (Quantavane).}
\end{corollary}

\begin{proof}
If $\varphi \colon P \to Q$ is an isomorphism of posets, then the map $\sigma \colon X_P \to X_Q$, given by $\sigma (x) := \varphi (x)$ and $\sigma (y, \, \varepsilon) := (\varphi(y), \, \varepsilon)$, is an isomorphism of solutions, by the same computation as in the proof of \thref{autDP}. Conversely, let $\sigma \colon X_P \to X_Q$ be an isomorphism of solutions. For $u \neq v$ in $X_P$ we have $D_P (u, \, v) \neq (u, \, v)$ if and only if one of the elements $u$, $v$ belongs to $P$ and the other one to $Z_P$; indeed, $D_P$ is the identity on $P \times P$ and on $Z_P \times Z_P$, while, for $x\in P$ and $z \in Z_P$, the first component of $D_P (x, \, z) = (d_x (z), \, x)$ belongs to $Z_P$ and the first component of $D_P (z, \, x)$ is $x$. The same holds for $D_Q$; hence $\sigma$ maps the partition $\{P, \, Z_P\}$ of $X_P$ onto the partition $\{Q, \, Z_Q\}$ of $X_Q$, and, since $|P| < |Z_P|$, we obtain $\sigma (P) = Q$ and $\sigma (Z_P) = Z_Q$. Let $\alpha := \sigma|_P \colon P \to Q$ and $\beta := \sigma|_{Z_P} \colon Z_P \to Z_Q$. Comparing the two sides of $(\sigma \times \sigma) \circ D_P = D_Q \circ (\sigma \times \sigma)$ on the pairs $(x, \, z) \in P \times Z_P$, we obtain $\beta \circ d_x = d_{\alpha(x)} \circ \beta$, for all $x\in P$. Now the argument in the proof of \thref{autDP}, with $Q$ in place of $P$ in the target, gives a bijection $\bar{\beta} \colon P \to Q$ such that $\bar{\beta} ({\downarrow} x) = {\downarrow} \alpha (x)$, for all $x \in P$, and, by \leref{downsets}, $\alpha = \bar{\beta}$ is an isomorphism of posets.
\end{proof}

The normal subgroup $(\bZ/2)^P$ of ${\rm Aut} (X_P, \, D_P)$ is due to the fact that the two elements $(y, \, 0)$ and $(y, \, 1)$ of $Z_P$ play symmetric roles, for each $y \in P$. In order to remove it, we use the involutive solutions associated to graphs in \cite[Theorem 2.14]{MilAut}.

Let $P$ be a poset and $B(P)$ the bipartite graph with set of vertices $Z_P = P \times \bZ/2$, in which $(y, \, 0)$ and $(x, \, 1)$ are adjacent if and only if $y \leq x$, and there are no other edges. The map $(y, \, x) \mapsto \{(y, \, 0), \, (x, \, 1)\}$ is a bijection between the set
\begin{equation*}
L_P := \{ (y, \, x) \in P \times P \, | \, y \leq x \}
\end{equation*}
and the set of edges of $B(P)$, and we identify these two sets. For $l = (y, \, x) \in L_P$ let $t_l \in {\rm Sym} (Z_P)$ be the transposition which interchanges $(y, \, 0)$ and $(x, \, 1)$. We denote by $E_P$ the involutive solution $R''_{B(P)}$ of the braid equation associated to the graph $B(P)$ in \cite[Theorem 2.14]{MilAut}, i.e. the bipartite solution associated to the datum $(L_P, \, Z_P, \, t)$, $t = (t_l)_{l\in L_P}$, on the set $X'_P := L_P \sqcup Z_P$. If $P$ is finite, then $|X'_P| = 3 |P| + m$, where $m$ is the number of pairs $(y, \, x) \in P \times P$ with $y < x$. Explicitly, for $l \in L_P$ and $z\in Z_P$ we have
\begin{equation}\eqlabel{EP}
E_P (l, \, z) = \bigl( t_l (z), \, l \bigr), \qquad E_P (z, \, l) = \bigl( l, \, t_l (z) \bigr),
\end{equation}
while $E_P$ is the identity on $L_P \times L_P$ and on $Z_P \times Z_P$. If $P \neq \emptyset$, then $B(P)$ has at least one edge, and by \cite[Theorem 2.14]{MilAut} the map
\begin{equation}\eqlabel{autBP}
{\rm Aut} \bigl( B(P) \bigr) \to {\rm Aut} (X'_P, \, E_P), \qquad \beta \mapsto \widehat{\beta},
\end{equation}
is an isomorphism of groups, where $\widehat{\beta}$ is equal to $\beta$ on $Z_P$, and $\widehat{\beta} (l)$ is the edge $\{\beta (y, \, 0), \, \beta (x, \, 1)\}$, for $l = (y, \, x) \in L_P$.

Recall that an \emph{anti-automorphism} of $P$ is a bijection $\varphi \colon P \to P$ such that $y \leq x$ if and only if $\varphi(x) \leq \varphi(y)$, for all $x$, $y\in P$, i.e. an isomorphism of posets $P \to P^{\rm op}$. We denote by ${\rm Aut}^{\pm} (P)$ the group of all pairs $(\varphi, \, s) \in {\rm Sym} (P) \times \bZ/2$ such that $\varphi$ is an automorphism of $P$ if $s = 0$ and an anti-automorphism of $P$ if $s = 1$, with the multiplication $(\varphi, \, s) (\varphi', \, s') := (\varphi \varphi', \, s + s')$. Then ${\rm Aut} (P) \cong \{ (\varphi, \, 0) \, | \, \varphi \in {\rm Aut} (P) \}$ is a subgroup of index $1$ or $2$ of ${\rm Aut}^{\pm} (P)$, and the index is $2$ if and only if $P \cong P^{\rm op}$. If $P$ has two comparable elements, then no bijection of $P$ is both an automorphism and an anti-automorphism, and ${\rm Aut}^{\pm} (P)$ is isomorphic to the group of all automorphisms and anti-automorphisms of $P$. Recall also that $P$ is \emph{connected} if its comparability graph is connected, i.e. $P$ is not the disjoint union of two non-empty subposets with no relations between them.

For $(\varphi, \, s) \in {\rm Aut}^{\pm} (P)$ we define the map $\sigma_{\varphi, s} \colon X'_P \to X'_P$ by
\begin{equation}\eqlabel{sigmaphis}
\sigma_{\varphi, s} (y, \, x) := \begin{cases} \bigl( \varphi(y), \, \varphi(x) \bigr), & \text{if} \,\,\, s = 0,\\ \bigl( \varphi(x), \, \varphi(y) \bigr), & \text{if} \,\,\, s = 1, \end{cases} \qquad \sigma_{\varphi, s} (u, \, \varepsilon) := \bigl( \varphi(u), \, \varepsilon + s \bigr),
\end{equation}
for all $(y, \, x) \in L_P$ and $(u, \, \varepsilon) \in Z_P$. Note that $\sigma_{\varphi, s}$ maps $L_P$ into $L_P$: if $y \leq x$, then $\varphi(y) \leq \varphi(x)$ for $s = 0$, and $\varphi(x) \leq \varphi(y)$ for $s = 1$. Moreover, $\sigma_{\varphi, s}$ is a bijection, whose inverse is $\sigma_{\varphi^{-1}, s}$.

\begin{theorem}\thlabel{autEP}
Let $P$ be a finite non-empty connected poset and $E_P$ the involutive solution of the braid equation on $X'_P = L_P \sqcup Z_P$ given by \equref{EP}. Then the map
\begin{equation*}
{\rm Aut}^{\pm} (P) \to {\rm Aut} (X'_P, \, E_P), \qquad (\varphi, \, s) \mapsto \sigma_{\varphi, s},
\end{equation*}
is an isomorphism of groups, where $\sigma_{\varphi, s}$ is given by \equref{sigmaphis}. In particular, if $P \not\cong P^{\rm op}$, then ${\rm Aut} (X'_P, \, E_P) \cong {\rm Aut} (P)$.
\end{theorem}

\begin{proof}
Since $P$ is non-empty, the graph $B(P)$ has at least one edge, and therefore the map \equref{autBP}, $\beta \mapsto \widehat{\beta}$, is an isomorphism of groups ${\rm Aut} (B(P)) \to {\rm Aut} (X'_P, \, E_P)$. For $(\varphi, \, s) \in {\rm Aut}^{\pm} (P)$ we denote by $\beta_{\varphi, s}$ the restriction of $\sigma_{\varphi, s}$ to $Z_P$, i.e.
\begin{equation*}
\beta_{\varphi, s} (u, \, \varepsilon) = \bigl( \varphi (u), \, \varepsilon + s \bigr), \qquad \text{for all} \,\, (u, \, \varepsilon) \in Z_P.
\end{equation*}
We shall prove that $(\varphi, \, s) \mapsto \beta_{\varphi, s}$ is an isomorphism of groups ${\rm Aut}^{\pm} (P) \to {\rm Aut} (B(P))$ and that $\widehat{\beta_{\varphi, s}} = \sigma_{\varphi, s}$; composing it with \equref{autBP}, we then obtain the first statement of the theorem.

We first show that $(\varphi, \, s) \mapsto \beta_{\varphi, s}$ is an injective morphism of groups ${\rm Aut}^{\pm} (P) \to {\rm Sym} (Z_P)$. For $(\varphi, \, s)$, $(\varphi', \, s') \in {\rm Aut}^{\pm} (P)$ and $(u, \, \varepsilon) \in Z_P$ we have
\begin{equation*}
\beta_{\varphi, s} \bigl( \beta_{\varphi', s'} (u, \, \varepsilon) \bigr) = \beta_{\varphi, s} \bigl( \varphi'(u), \, \varepsilon + s' \bigr) = \bigl( \varphi \varphi'(u), \, \varepsilon + s' + s \bigr) = \beta_{\varphi \varphi', \, s + s'} (u, \, \varepsilon),
\end{equation*}
i.e. $\beta_{\varphi, s} \circ \beta_{\varphi', s'} = \beta_{(\varphi, s) (\varphi', s')}$. In particular, $\beta_{\varphi, s}$ is a bijection, with inverse $\beta_{\varphi^{-1}, s}$. Moreover, $(\varphi, \, s)$ is determined by $\beta_{\varphi, s}$, since $\beta_{\varphi, s} (u, \, 0) = (\varphi(u), \, s)$, for all $u\in P$.

Next, we show that $\beta_{\varphi, s} \in {\rm Aut} (B(P))$ and $\widehat{\beta_{\varphi, s}} = \sigma_{\varphi, s}$. Let $(y, \, x) \in L_P$, i.e. $y \leq x$, and consider the corresponding edge $\{(y, \, 0), \, (x, \, 1)\}$ of $B(P)$. If $s = 0$, then $\varphi$ is an automorphism of $P$, hence $\varphi(y) \leq \varphi(x)$, and $\beta_{\varphi, 0}$ maps our edge to
\begin{equation*}
\bigl\{ (\varphi(y), \, 0), \, (\varphi(x), \, 1) \bigr\},
\end{equation*}
which is the edge $(\varphi(y), \, \varphi(x)) = \sigma_{\varphi, 0} (y, \, x)$. If $s = 1$, then $\varphi$ is an anti-automorphism of $P$, hence $\varphi(x) \leq \varphi(y)$, and $\beta_{\varphi, 1}$ maps our edge to
\begin{equation*}
\bigl\{ (\varphi(y), \, 1), \, (\varphi(x), \, 0) \bigr\} = \bigl\{ (\varphi(x), \, 0), \, (\varphi(y), \, 1) \bigr\},
\end{equation*}
which is the edge $(\varphi(x), \, \varphi(y)) = \sigma_{\varphi, 1} (y, \, x)$. Thus $\beta_{\varphi, s}$ maps edges to edges. The same holds for its inverse $\beta_{\varphi^{-1}, s}$, since $(\varphi^{-1}, \, s) \in {\rm Aut}^{\pm} (P)$; hence $\beta_{\varphi, s}$ is an automorphism of $B(P)$. Moreover, $\widehat{\beta_{\varphi, s}}$ and $\sigma_{\varphi, s}$ are both equal to $\beta_{\varphi, s}$ on $Z_P$, and the above computation shows that they coincide on $L_P$; thus $\widehat{\beta_{\varphi, s}} = \sigma_{\varphi, s}$.

It remains to show that any automorphism of $B(P)$ is of the form $\beta_{\varphi, s}$. We first show that the graph $B(P)$ is connected. For any $u\in P$, the vertices $(u, \, 0)$ and $(u, \, 1)$ are adjacent, as $u \leq u$. If $u$, $v \in P$ are comparable, say $u \leq v$, then $(u, \, 0)$ and $(v, \, 1)$ are adjacent, and hence the four vertices $(u, \, 0)$, $(u, \, 1)$, $(v, \, 0)$, $(v, \, 1)$ belong to the same connected component of $B(P)$. Since $P$ is connected, any two elements $u$, $v\in P$ can be joined by a sequence $u = u_0, \, u_1, \, \dots, \, u_k = v$ of elements of $P$ such that $u_{i-1}$ and $u_i$ are comparable, for $i = 1, \, \dots, \, k$; thus all the vertices of $B(P)$ belong to the same connected component.

Now let $\beta \in {\rm Aut} (B(P))$. The graph $B(P)$ is bipartite, with parts $P \times \{0\}$ and $P \times \{1\}$. Since $\beta$ maps edges to edges, $\beta (P \times \{0\})$ and $\beta (P \times \{1\})$ give again a decomposition of $B(P)$ as a bipartite graph; as $B(P)$ is connected, this decomposition is unique up to the order of the two parts (see \seref{prel}). Hence $\beta$ either preserves both parts or interchanges them, and we discuss the two cases separately.

Assume first that $\beta$ preserves the two parts. Then there exist two bijections $\gamma$, $\delta \colon P \to P$ such that
\begin{equation*}
\beta (u, \, 0) = \bigl( \gamma (u), \, 0 \bigr) \qquad {\rm and} \qquad \beta (u, \, 1) = \bigl( \delta (u), \, 1 \bigr), \qquad \text{for all} \,\, u\in P.
\end{equation*}
As $\beta$ is an automorphism of $B(P)$, $\{(y, \, 0), \, (x, \, 1)\}$ is an edge if and only if $\{(\gamma(y), \, 0), \, (\delta(x), \, 1)\}$ is an edge, i.e. $y \leq x$ if and only if $\gamma (y) \leq \delta (x)$, for all $x$, $y\in P$. This means that $\gamma ({\downarrow} x) = {\downarrow} \delta (x)$, for all $x\in P$. By \leref{downsets}, applied to the bijections $\gamma$ and $\delta$, we obtain that $\gamma = \delta =: \varphi$ is an automorphism of $P$, and then $\beta = \beta_{\varphi, 0}$.

Assume now that $\beta$ interchanges the two parts. Then there exist two bijections $\gamma$, $\delta \colon P \to P$ such that
\begin{equation*}
\beta (u, \, 0) = \bigl( \gamma (u), \, 1 \bigr) \qquad {\rm and} \qquad \beta (u, \, 1) = \bigl( \delta (u), \, 0 \bigr), \qquad \text{for all} \,\, u\in P.
\end{equation*}
As above, $\{(y, \, 0), \, (x, \, 1)\}$ is an edge if and only if $\{(\delta(x), \, 0), \, (\gamma(y), \, 1)\}$ is an edge, i.e. $y \leq x$ if and only if $\delta (x) \leq \gamma (y)$, for all $x$, $y\in P$. For $y\in P$, the down-set of $y$ in $P^{\rm op}$ is $\{ x\in P \, | \, y \leq x \}$; hence the last equivalence means that $\delta$ maps the down-set of $y$ in $P^{\rm op}$ onto the down-set ${\downarrow} \gamma (y)$ of $\gamma(y)$ in $P$, for all $y\in P$. By \leref{downsets}, applied to the posets $P^{\rm op}$ and $P$ and to the bijections $\delta$ and $\gamma$, we obtain that $\delta = \gamma =: \varphi$ is an isomorphism of posets $P^{\rm op} \to P$, i.e. an anti-automorphism of $P$. Then $(\varphi, \, 1) \in {\rm Aut}^{\pm} (P)$ and $\beta = \beta_{\varphi, 1}$.

Therefore, $(\varphi, \, s) \mapsto \beta_{\varphi, s}$ is an isomorphism of groups ${\rm Aut}^{\pm} (P) \to {\rm Aut} (B(P))$, and, composing it with the isomorphism \equref{autBP}, we obtain that $(\varphi, \, s) \mapsto \widehat{\beta_{\varphi, s}} = \sigma_{\varphi, s}$ is an isomorphism of groups ${\rm Aut}^{\pm} (P) \to {\rm Aut} (X'_P, \, E_P)$. Finally, if $P \not\cong P^{\rm op}$, then $P$ has no anti-automorphisms, hence ${\rm Aut}^{\pm} (P) = \{ (\varphi, \, 0) \, | \, \varphi \in {\rm Aut} (P) \} \cong {\rm Aut} (P)$, and ${\rm Aut} (X'_P, \, E_P) \cong {\rm Aut} (P)$.
\end{proof}

\begin{remark}\relabel{invol}
The connectedness of $P$ cannot be dropped in \thref{autEP}: if $P$ is an antichain with $n \geq 2$ elements, then $B(P)$ is the disjoint union of $n$ edges, and ${\rm Aut} (X'_P, \, E_P) \cong {\rm Aut} (B(P)) \cong (\bZ/2)^n \rtimes S_n$ has order $2^n \, n!$, while ${\rm Aut}^{\pm} (P)$ has order $2 \, n!$.

On the contrary, the solutions $E_P$ do not distinguish a poset from its opposite: for any finite poset $P$, the map $X'_P \to X'_{P^{\rm op}}$ given by $(u, \, \varepsilon) \mapsto (u, \, \varepsilon + 1)$ on $Z_P$ and $(y, \, x) \mapsto (x, \, y)$ on $L_P$ is an isomorphism of solutions $(X'_P, \, E_P) \cong (X'_{P^{\rm op}}, \, E_{P^{\rm op}})$, induced by the isomorphism of graphs $B(P) \cong B(P^{\rm op})$ which interchanges the two parts.
\end{remark}

\begin{example}\exlabel{C3}
We illustrate the results of this paper on the poset $P := \{ a_i, \, b_i, \, c_i \, | \, i \in \bZ/3 \}$ with nine elements, ordered by
\begin{equation}\eqlabel{exC3}
a_i < b_i, \qquad b_i < c_i, \qquad a_i < c_i, \qquad a_i < c_{i+1}, \qquad \text{for all} \,\, i \in \bZ/3,
\end{equation}
and with no other strict relations (the indices are taken modulo $3$). This poset is used in the proof of \cite[Proposition 2.3]{BarmakBarreto24}, which shows that $9$ is the smallest number of elements of a poset with cyclic automorphism group of order $3$. Its automorphism group is
\begin{equation*}
{\rm Aut} (P) = \{ {\rm Id}_P, \, \rho, \, \rho^2 \} \cong \bZ/3,
\end{equation*}
where $\rho \colon P \to P$ is given by $\rho (a_i) := a_{i+1}$, $\rho (b_i) := b_{i+1}$ and $\rho (c_i) := c_{i+1}$.

On the other hand, the map $\psi \colon P \to P$ given by $\psi (a_i) := c_{-i}$, $\psi (b_i) := b_{-i}$ and $\psi (c_i) := a_{-i}$ is an anti-automorphism of $P$: it reverses each of the relations \equref{exC3}; for instance, $a_i < c_{i+1}$ is mapped to $a_{-i-1} < c_{-i}$. Thus $P \cong P^{\rm op}$, the anti-automorphisms of $P$ are $\psi \rho^k$, $k = 0$, $1$, $2$, and, since $\psi \rho \psi^{-1} = \rho^{-1}$, the
group ${\rm Aut}^{\pm} (P)$ is isomorphic to the symmetric group $S_3$.

The solution $R_P$ lives on a set with $9$ elements: $R_P (x, \, y) = (y, \, x)$ if $(x, \, y)$ is one of the $12$ pairs
\begin{equation*}
(b_i, \, a_i), \qquad (c_i, \, b_i), \qquad (c_i, \, a_i), \qquad (c_{i+1}, \, a_i), \qquad i \in \bZ/3,
\end{equation*}
and $R_P (x, \, y) = (x, \, y)$ for the other $69$ pairs. By \thref{RP}, $R_P$ is an idempotent solution of the braid equation and ${\rm Aut} (P, \, R_P) = {\rm Aut} (P) \cong \bZ/3$.

The involutive solution $D_P$ lives on the set $X_P = P \sqcup (P \times \bZ/2)$ with $27$ elements and, by \thref{autDP}, ${\rm Aut} (X_P, \, D_P) \cong (\bZ/2)^9 \rtimes \bZ/3$ is a group of
order $2^9 \cdot 3 = 1536$.

The poset $P$ is connected: $a_i$ is comparable to $b_i$, $c_i$ and $c_{i+1}$, for all $i\in \bZ/3$. The set $L_P$ has $9 + 12 = 21$ elements, hence the involutive solution $E_P$ lives on the set $X'_P = L_P \sqcup Z_P$ with $21 + 18 = 39$ elements. By \thref{autEP} and the above description of ${\rm Aut}^{\pm} (P)$, we have
\begin{equation*}
{\rm Aut} (X'_P, \, E_P) \cong {\rm Aut}^{\pm} (P) \cong S_3,
\end{equation*}
the group being generated by $\sigma_{\rho, 0}$ and $\sigma_{\psi, 1}$, where, by \equref{sigmaphis}, $\sigma_{\psi, 1} (u, \, \varepsilon) = (\psi(u), \, \varepsilon + 1)$ and $\sigma_{\psi, 1} (y, \, x) = (\psi(x), \, \psi(y))$. Thus, in this case, the automorphism group of $E_P$ is larger than ${\rm Aut} (P)$: the anti-automorphism $\psi$ of $P$ produces automorphisms of $E_P$ which interchange $P \times \{0\}$ and $P \times \{1\}$.

In order to obtain an involutive solution whose automorphism group is exactly $\bZ/3$, we add a greatest element to $P$: let $P' := P \sqcup \{ \top \}$, with $x < \top$ for all $x\in P$. Any automorphism of $P'$ fixes $\top$, hence ${\rm Aut} (P') \cong {\rm Aut} (P) \cong \bZ/3$. The poset $P'$ is connected, and $P' \not\cong P'^{\rm op}$, since $P'$ has a greatest element but no least element (its minimal elements are $a_0$, $a_1$ and $a_2$). As $P'$ has $12 + 9 = 21$ pairs $(y, \, x)$ with $y < x$, the set $X'_{P'}$ has $3 \cdot 10 + 21 = 51$ elements, and by \thref{autEP} we obtain
\begin{equation*}
{\rm Aut} (X'_{P'}, \, E_{P'}) \cong {\rm Aut} (P') \cong \bZ/3.
\end{equation*}
\end{example}

%%%%%%%%%%%%%%%%%%%%%%%%%%%%%%%%%%%%%%%%%%%%%%%%%%%%%%%%%%%%%%%%%%%%%%%%%%%%%

\subsection*{Acknowledgements}
The author thanks Arhan Canli (quantavane.com) for an independent computer verification of \thref{RP}, \thref{flipfix} and \reref{flip}, and for the question answered by \coref{isoDP}. Claude Opus~5.5 (Anthropic) was used throughout the preparation of this paper to find references for the problem under consideration, improving the author’s English, to explore constructions and proofs, proofreading, correcting typos, and making improvements. All statements, proofs and computational results were independently checked by the author, who takes sole responsibility for them.

\subsection*{Conflict of interest statement}
The author declares no conflict of interest.

\subsection*{Data availability statement}
Data sharing is not applicable to this article as no datasets were generated or analysed during the current study.

\newcommand{\Addresses}{{%
  \bigskip
  \footnotesize
  \textsc{Faculty of Mathematics and Computer Science, University of Bucharest}
  \par\nopagebreak
  \textsc{Str. Academiei 14, RO-010014 Bucharest 1, Romania}
  \par\nopagebreak
  \textit{E-mail address}: \texttt{gigel.militaru@gmail.com}
}}

\Addresses

\end{document}